\documentclass[12pt]{article}
\usepackage{amsmath,amssymb,amsthm}
\usepackage{graphicx}
\usepackage[T1]{fontenc}

\newtheorem{thm}{\bf Theorem}

\newtheorem{rem}{\bf Remark}
\newtheorem{definition}[thm]{\bf Definition}
\newtheorem{hypothesis}[thm]{\bf Hypothesis}

\theoremstyle{definition}
\newtheorem{ex}{\bf Example}

\begin{document}

\title{Special cases of pairwise comparisons \\ matrices represented by Toeplitz matrices}

\author{Viera \v{C}er\v{n}anov\'a \footnote {Slovak University of Technology, Institute of Computer Science and Mathematics, Ilkovi\v{c}ova 3, 812 19 Bratislava, Slovak Republic, vieracernanova@hotmail.com} and Waldemar W. Koczkodaj \footnote{Computer Science, Laurentian University, Sudbury, Canada, {wkoczkodaj@cs.laurentian.ca}}}

\date{March, 2017}

\maketitle
\begin{abstract}
This study presents special cases of inconsistent pairwise comparisons PC matrices and analysis of their eigenvalue-based inconsistency index using mathematical methods.
	All studied special cases of PC matrices are Toeplitz matrices with only three different entries $1$, $x$, and $1/x$. A new type of circulant pairwise comparisons matrix has been introduced.
	Although this class of PC matrices may be perceived as restricted, it is general enough to cover numerous levels of eigenvalue-based inconsistency index
	from the lowest to the highest. Both exact mathematical expressions and estimations, where the exact expression was impossible to find, are provided. \\


\end{abstract}

\section{INTRODUCTION}
\label{intro}
In science and common practice, we compare entities in pairs, often without
realizing it. For this reason, pairwise comparisons might have been one of the first scientific methods. The first application of pairwise comparisons to elections has been recently attributed to the works of Ramon Llull, the 12-th century  mathematician, logician, philosopher, Majorcan writer, and mystic \cite{aboutLlull}. However, it is easy to envision the practical use of pairwise comparisons taking place in the Stone Age. Stones must have been compared to each other in pairs to decide which is better suited to use as a tool.

When we have no unit (e.g., for reliability), we may consider the construction of a pairwise comparisons  matrix (PC matrix) to express our assessments based on relative comparisons of  attributes (such as safety or reliability). After all, common sense and an old adage (commonly attributed to Creighton Abrams) calls for ``take one bite at a time'' when it comes to eating an elephant. Our elephant is the processing of subjective data, especially for the decision making where the ``satisfying'' approach is often used.  
Herbert A. Simon, the Nobel prize winner, proposed bounded rationality (``satisfying'') 
as a vital alternative to the exclusiveness of using mathematical theory for decision making. Pairwise comparisons supports the concept of bounded rationality.

In this study, we consider special cases of pairwise comparisons matrices (abbreviated to PC matrix here), which are $n\times n$ reciprocal matrices $\mathbf{A}=\left( {{a}_{ij}} \right)$ with positive entries. PC matrix $\mathbf{A}$ is called \emph{reciprocal} if ${{a}_{ij}}=1/{{a}_{ji}}$ for $i,j=1,...,n$. Obviously, ${{a}_{ii}}=1$ for any $i$, and thus the trace is $Tr\left( \mathbf{A} \right)=n.$ However, blind wine testing may result not only in the lack of reciprocity but even to the lack of $1$s on the main diagonal since comparing the same wine to itself (especially at the end of the tasting day), may not necessarily be correct. Due to the Internet, different assessors may compare projects or at least their components in different locations. In such situations, it is even anticipated that some (if not most) assessments may not be reciprocal.

A pairwise comparisons matrix $\mathbf{A}$ is called \emph{consistent} (or \emph{transitive}) if 
\begin{equation} \label{eq:cc} 
{{a}_{ij}}\cdot {{a}_{jk}}={{a}_{ik}} \text{ for } i,j,k=1,...,n. 
\end{equation} 
We refer to this as the  consistency condition. While every consistent matrix is reciprocal, the converse is false in general. If the consistency condition does not hold, the matrix is inconsistent (or intransitive). Given a reciprocal $n\times n$ matrix $\mathbf{A}$ which is not consistent, the theory attempts to provide a consistent $n\times n$ matrix $\mathbf{A}'$ which differs from $\mathbf{A}$ ``as little as possible''. The challenge of the pairwise comparisons method comes from the lack of consistency of pairwise comparisons matrices which arise in practice.

The quotients ${{a}_{ij}}={{E}_{i}}/{{E}_{j}}$ express a relative preference of an entity $E_{i}$ over $E_{j}$. An entity could be any object, attribute of it or a stimulus. Consistent matrices correspond to the ideal situation in which there are exact values ${{E}_{1}},...,{{E}_{n}}$ for the stimuli, since ${{a}_{ij}}={{E}_{i}}/{{E}_{j}}$ form a consistent matrix for all (even random) positives values ${{E}_{i}}$. It is an important observation since the implication of it is that, a problem of approximation is really a~problem of selecting a norm and a distance minimization. Notice that the vector $\left( {{E}_{1}},...,{{E}_{n}} \right)$ is unique up to a multiplicative constant. For the Euclidean norm, the vector of geometric means (which is equal to the principal eigenvector for a consistent PC matrix) is the one which generates it. The seminal
study \cite{Saaty1977} had a profound impact on the pairwise
comparisons research. It has strongly endorsed the use of the eigenvector,
corresponding to the principal eigenvalue, for approximation of a~given inconsistent but reciprocal PC matrix. Numerous studies show the lack of evidence for the superiority of the eigenvector solution. It is expressed in the highly cited \cite{BV2008} and a sizable collaboration (recently published) \cite{collab}. 

In layman's terms, when we have three entities: $A$,$B$,$C$, then $(A/B)\times (B/C)$ must yield the same result as the comparison $A/C$. However,  more often than not, relation~\eqref{eq:cc} does not hold when these three comparisons are carried independently. When the comparisons (e.g., important, safer, etc.) or entities (safety, reliability, etc.) are subjective, the inconsistency is unavoidable in practice. In fact, the lack of inconsistency in such cases may be suspicious but it does not mean that inconsistency is desirable or should be tolerated. The ``GIGO rule'' (GIGO stands ``garbage-in, garbage-out'') in the field of computer science and information technology refers to the fact that output quality depends on the quality of input data.

At least three studies: \cite{GS1958,KB1939,S1961} defined and examined inconsistency in pairwise comparisons between 1939 and 1961. The importance of inconsistency in pairwise comparisons is undeniable since it is still studied. Pairwise comparisons are used in the decision making process for choosing from two alternatives (yes/no, left/right, increase/decrease...). In particular, reciprocity and consistency is studied in \cite{RV2013}. 

Toeplitz matrices occur in a variety of problems in applied mathematics, engineering theory or computer science, such as integral equations, time series analysis, signal processing, theory of codes, cryptography etc. They have been examined  for a long time, both in theory and application. Nevertheless, they seem to offer an inexhaustible source of research ideas, and domains of application. Recently, it was shown that every square matrix over $\mathbb{C}$ is a~product of Toeplitz matrices \cite{YL2016}. We read in \cite{YL2016}:
\begin{quote}
The fact that a matrix may be expressed as a product of a lower-triangular with an upper-triangular matrix (LU), or of an orthogonal with an upper-triangular matrix (QR), or of two orthogonal matrices with diagonal one (SVD) is a cornerstone of modern numerical computations.
\end{quote}
The decomposition of a matrix in a product of Toeplitz matrices provides an alternative to (LU), (QR) or (SVD) decomposition. 

In this work, we examine levels of inconsistency in some special cases of PC matrices from the class of Toeplitz matrices. 

The  paper is organized as follows. In Section \ref{sec:pre}, some definitions and properties, relating to used matrix methods are outlined together with some basic concepts of pairwise comparisons.
Sections \ref{CC} to \ref{sec:FPC} are devoted to various special types of pairwise comparisons matrices, which are $n\times n$ Toeplitz matrices with entries $1$, $x$, $1/x$. In most cases a closed form of the consistency index \eqref{eq:CI}  is found as a function of $x$ and $n$. An open problem formulated in subsection~\ref{LPCeven} as Hypothesis~\ref{hypo} proposes an estimate of the consistency indicator of layer-cake PC matrices with even order.

\section{Preliminaries}
\label{sec:pre}
Let us recall from \cite{KS2014} that a vector $\left( {{w}_{1}},...,{{w}_{n}} \right)$, with even random positive real coordinates, generates a consistent PC matrix $\mathbf{A}=\left( {{a}_{ij}} \right)$ with ${{a}_{ij}}={{w}_{i}}/{{w}_{j}}$ for $i,j=1,...,n$.

The following observations are known properties of consistent matrices (some of them are introduced in \cite{Saaty1977}), and trivial to check:

\begin{rem}
\label{PCpropreties} \
\begin{enumerate}
\renewcommand\labelenumi{(\roman{enumi})}  
%\item Arbitrary consistent matrix is reciprocal. The inverse is not true for $n\ge 3$, but any reciprocal  $2\times 2$ matrix is consistent.
\item Each row of a consistent matrix is a constant multiple of any other row. Similarly, each column of a consistent matrix is a constant multiple of any other column.
\item $0$ is an eigenvalue with multiplicity $n-1$ of any consistent $n\times n$ matrix, 
hence the unique non-zero eigenvalue of a consistent $n\times n$ matrix is $n.$
\item Every column of a consistent PC matrix is its eigenvector corresponding to the eigenvalue $n.$
\end{enumerate}
\end{rem}

In \cite{AG1993}, the greatest lower bound and the least upper bound were found for the Perron root (called principal or Perron's eigenvalue, and equal to the spectral radius of the matrix with non-negative entries) of a PC matrix:
\begin{quote}
Theorem [quoted ``as is'' from \cite{AG1993}]: Let $\mathbf{A}$ be a positive reciprocal matrix with entries $1/x\le {{a}_{ij}}\le x,1\le i,j\le n$ for some $x\ge 1$, and let ${{\lambda }_{\max }}$ denote the largest eigenvalue of $\mathbf{A}$ in modulus, which is known to be real and positive from the Perron-Frobenius theorem. Then \[n\le {{\lambda }_{\max }}\le 1+\tfrac{1}{2}\left( n-1 \right)\left( x+\tfrac{1}{x} \right),\] the lower and upper bound being reached if and only if $\mathbf{A}$ is supertransitive or maximally intransitive, respectively.
\end{quote}

Classical bounds for the Perron root of a nonnegative matrix are the minimum and maximum row sums (see, for example Theorem 8.1.22 in \cite{HJ1985}):
\begin{quote}
%\begin{thm}
\label{HornJoh}
%Perron root ${{\lambda }_{\max }} $ of a non-negative matrix 
%$\mathbf{A}={{({{a}_{ij}})}_{n\times n}}$ satisfies 
\[\underset{i}{\mathop{\min }}\,\sum\nolimits_{j=1}^{n}{{{a}_{ij}}}\le {{\lambda }_{\max }}\le \underset{i}{\mathop{\max }}\,\sum\nolimits_{j=1}^{n}{{{a}_{ij}}}.\]
%\end{thm}
\end{quote}

According to remark~\ref{PCpropreties}{ (ii)}, $n$ is the Perron root of any consistent $n\times n$ PC matrix. The scaled (by $n-1$) difference between the principal eigenvalue ${{\lambda }_{\max }}$ of a given reciprocal matrix $\mathbf{A}$ and that of a consistent matrix of the same order $n$ was defined in \cite{Saaty1977} and erroneously named  \emph{consistency indicator} $CI\left( \mathbf{A} \right)$:
\begin{equation} \label{eq:CI}
CI\left( \mathbf{A} \right)=\frac{{{\lambda }_{\max }}-n}{n-1}
\end{equation}
since it leads to error tolerance problems. These problems were demonstrated in \cite{KS2014} by two counter-examples of PC matrices CPC and FPC type (here discussed). As we will show, the scaling by $n-1$ in \cite{Saaty1977} is highly unfortunate since a better scaling exists and could be used for normalizing $CI\left( \mathbf{A} \right)$ to $[0,1]$.

In accordance with the above cited Theorem in \cite{AG1993}, and with the use of \eqref{eq:CI} we get 
\begin{equation}\label{eq:CI_bounds}
0\le CI\left( \mathbf{A} \right)\le \frac{1+\frac{1}{2}\left( n-1 \right)\left( x+\frac{1}{x} \right)-n}{n-1}=\frac{1}{2}\left( x+\frac{1}{x} \right)-1=\frac{{{\left( x-1 \right)}^{2}}}{2x}.
\end{equation}

According to \cite{AG1993}, a PC matrix $\mathbf{A}$ is maximally inconsistent (maximally intransitive) if and only if $CI\left( \mathbf{A} \right)=\tfrac{{{\left( x-1 \right)}^{2}}}{2x}$. \\

In \cite{K1993}, a \emph{distance-based inconsistency index (Kii)} was defined. It was generalized in \cite{DK1994} and simplified in \cite{KS2014} to:
\begin{equation}
Kii\left( \mathbf{A} \right)=\max_{i,j,k \le dim \left( \mathbf{A} \right)}\left(1-\min\left(\frac{a_{ik}}{a_{ij}a_{jk}},\frac{a_{ij}a_{jk}}{a_{ik}}\right)\right).
\end{equation}

Notice, that $CI$ index is nonnegative and may be arbitrarily large, while $Kii$ has values in $[0,1)$. Evidently, both are equal to $0$ if and only if the PC matrix is consistent (supertransitive). 

In \cite{KS2014}, two counter-examples, with the mathematical reasoning and proofs, were provided that the eigenvalue-based inconsistency tolerates an error of any arbitrarily large value (e.g., 1,000,000\% or whatever our imagination calls for). 

Recently, two major studies \cite{collab} and \cite{K2015paradox} have been published. They address problems related to the use of eigenvector-based solution and rating scales for data entry, without normalization, in PC matrices.

All considered matrices in this study are positive real Toeplitz matrices. The following Definition~\ref{Toeplitz} (due to Toeplitz) introduces them. %to \ref{def_gers}. % restricted to real case.
\begin{definition}
\label{Toeplitz}
A square matrix $\mathbf{A}=\left( {{a}_{ij}} \right)$ of order $n$ is called real finite Toeplitz matrix if ${{a}_{ij}}={{c}_{j-i}}$ for any real constants ${{c}_{1-n}},...,{{c}_{n-1}}.$
\end{definition}

A Toeplitz matrix is called also a diagonal-constant matrix since each descending diagonal from left to right is constant. 

\begin{definition}
A circulant matrix is a finite Toeplitz matrix of the following form:
$$
\begin{bmatrix}
{{c}_{0}} & {{c}_{1}} & {{c}_{2}} & \ldots  & {{c}_{n-1}}  \\
{{c}_{n-1}} & {{c}_{0}} & {{c}_{1}} & \cdots  & {{c}_{n-2}}  \\
{{c}_{n-2}} & {{c}_{n-1}} & {{c}_{0}} & \cdots  & {{c}_{n-3}}  \\
\vdots  & \vdots  & \vdots  & \ddots  & \vdots   \\
{{c}_{1}} & {{c}_{2}} & {{c}_{3}} & \cdots  & {{c}_{0}}  \\
\end{bmatrix}, $$
where each row is a cyclic shift of the row above it.
\end{definition}

It is known, that the eigenvalues of a circulant $n\times n$ matrix can be written
\begin{equation}\label{eq:eval_circulant}
{{\lambda }_{m}}=\sum\limits_{k=0}^{n-1}{{{c}_{k}}{{e}^{-2\pi imk/n}}},m=1,...,n. 	
\end{equation}

\begin{definition}
The Frobenius norm of a real finite $n\times n$ matrix $\mathbf{A}=\left( {{a}_{ij}} \right)$ with eigenvalues ${{\lambda }_{i}},i=1,..,n$ is
\begin{equation}\label{eq:Frobenius} 
{{\left\| \mathbf{A} \right\|}_{F}}={{\left( \sum\limits_{i,j=1}^{n}{{{\left| {{a}_{ij}} \right|}^{2}}} \right)}^{1/2}}={{\left( {Tr\left( {{\mathbf{A}}^{T}}\mathbf{A} \right)} \right)}^{1/2}}.
\end{equation}
\end{definition}

\begin{definition}[from \cite{G1931} ]
\label{def_gers}
Let $\mathbf{A}=\left( {{a}_{ij}} \right)$ be a real matrix. For $i=1,...,n$ denote ${{R}_{i}}=\sum\nolimits_{j\ne i}{\left| {{a}_{ij}} \right|}$. Let $D\left( {{a}_{ii}},{{R}_{i}} \right)$ be the closed disc centered at ${{a}_{ii}}$ with radius ${{R}_{i}}$. Such a disc is called a~Gerschgorin disc.
\end{definition}

\begin{thm}[from \cite{G1931}]
Every eigenvalue of $\mathbf{A}$ lies within at least one of the Gerschgorin discs $D\left( {{a}_{ii}},{{R}_{i}} \right)$.
\end{thm}

\section{A consistent PC matrix with all entries equal to 1 (CC)}
\label{CC}
$\mathbf{CC}\left( n \right)=\left( {{a}_{ij}} \right)$ with ${{a}_{ij}}=1$ for $i,j=1,...,n$ is consistent. $\mathbf{CC}\left( n \right)$ is the simplest consistent PC matrix, used often in  pairwise comparisons. According to the Remark~\ref{PCpropreties}, $0$ is an eigenvalue with multiplicity $n-1$, the Perron root is $n$ and the consistency index~$0$. The 1-dimensional eigenspace corresponding to the eigenvalue $n$ is spanned by the vector ${{\left( 1,1,...,1 \right)}^{T}}$. 

This matrix is an important case since it is usually assumed for the initial state: everything is ``equal'' (or unknown) each other. Changing even one value in this matrix may cause problems of considerable importance as the next section demonstrates.

\section{Corner PC matrix (CPC}
\label{CPC}
In this section, a reciprocal matrix with all $1$'s except for two terms ${{a}_{ij}}=x$ and ${{a}_{ji}}=1/x$ is analyzed. First, notice that all such $n\times n$ reciprocal PC matrices with equal $x$ have the same characteristic polynomial. This takes place since every PC matrix of this type is obtained by one or two row (and equal number of column) exchanges from another one, say $\mathbf{A}$, which preserves the determinant $\det \left( \mathbf{A}-\lambda \mathbf{I} \right)$. Consequently, all matrices from this family have the equal spectrum and consistency indicator. This property is not a surprise since the position of $x$ (and its reciprocal, $1/x$) in the matrix does not influence the number of triads for a given $n$.

Let us denote this family of matrices by ${{\mathcal{F}}_{o}}$.
As a representative of ${{\mathcal{F}}_{o}}$, let us choose a corner pairwise comparisons matrix: $\mathbf{CPC}\left( x,n \right)$ with $n\ge 3$ and $x>1$ defined by: \\
\[
\mathbf{CPC}\left( x,n \right)=\begin{bmatrix}
1 & 1 & \cdots  & 1 & x  \\
1 & 1 & \cdots  & 1 & 1  \\
\vdots  & \vdots  & \ddots  & \vdots  & \vdots   \\
1 & 1 & \cdots  & 1 & 1  \\
1/x & 1 & \cdots  & 1 & 1  \\
\end{bmatrix}.
\]

Evidently, there are $n-2$ equal rows, therefore $0$ is an eigenvalue with multiplicity $n-3$. The characteristic polynomial has a form: \\
\begin{equation}
\label{eq:charpol}
{{\left( -1 \right)}^{n-1}}{{\lambda }^{n-3}}\left( {{\lambda }^{3}}-n{{\lambda }^{2}}-\left( n-2 \right)\left( x-2+\tfrac{1}{x} \right) \right),
\end{equation}
thus the non-zero eigenvalues can be computed using Cardano's method. The matrix $\mathbf{CPC}\left( x,n \right)$ has one real positive and two complex conjugate non-zero eigenvalues.  The Perron root is
\begin{equation}
\label{eq:lambdaCPC}
{{\lambda }_{\max }}\left( x,n \right)=B\left( x,n \right)+\frac{n^2}{9B\left( x,n \right)}+\frac{n}{3},
\end{equation}
where $B\left( x,n \right)$ is
\[
{{\left( \frac{{{n}^{3}}}{27}+\frac{\left( x-1 \right)\sqrt{4{{n}^{3}}\left( n-2 \right)x+27{{\left( n-2 \right)}^{2}}{{\left( x-1 \right)}^{2}}}}{2\cdot 3\sqrt{3}x}+\frac{\left( n-2 \right){{\left( x-1 \right)}^{2}}}{2x} \right)}^{1/3}}.
\]
This expression for ${\lambda }_{\max }$ would be unusually complicated for our aim, which is an estimate of $CI$ yet it is simpler than in \cite{KS2014}. 

In \cite{KS2014}, a very good estimate of ${\lambda }_{\max }$ was produced as:  
\begin{equation}
\label{eq:1odhadCICPC}
CI\left(\mathbf{CPC}\left(x,n \right) \right) \le \frac{x}{n^2}.
\end{equation}
\noindent The estimate \eqref{eq:1odhadCICPC} may be further improved. Using \eqref{eq:charpol} and \cite{KS2014}, we obtain
\begin{equation}
\label{eq:2odhadCICPC}
CI\left(\mathbf{CPC}\left(x,n \right) \right) \le \frac{n-2}{{{n}^{2}}\left( n-1 \right)}\frac{{{\left( x-1 \right)}^{2}}}{x},
\end{equation}
\noindent which corresponds to the inequality (2) in \cite{KS2014}.\\
\noindent As expected, $CI$ of $\mathbf{CPC}\left( x,n \right)$ matrix is converging to 0 for $n \rightarrow \infty$ when $x$ is fixed (no matter to how large value). In the mathematical terminology terms, it satisfies: $$\underset{n\to \infty }{\mathop{\lim }}\,CI\left(\mathbf{CPC}\left( x,n \right) \right)=0.$$

\section{Layer-cake PC matrix (LPC)}
\label{sec:LPC}
Let us shift our attention from PC matrices with the smallest Perron root to PC matrices with the largest Perron root. For this aim, let us consider $x>0$, $x\ne 1$ and a~reciprocal $n\times n$ matrix $\mathbf{LPC}\left( x,n \right)$ with diagonal elements ${{a}_{ii}}=1$ and entries $x$ and $1/x$ alternating on diagonals located above the main diagonal as illustrated by the following layer-cake PC matrix:
\[
\mathbf{LPC}\left( x,n \right)=\begin{bmatrix}
1 & x & 1/x & x & \cdots   \\
1/x & 1 & x & 1/x & \cdots   \\
x & 1/x & 1 & x & \cdots   \\
\vdots  & \vdots  & \vdots  & \ddots  & \cdots   \\
\end{bmatrix}.\]

The number of $1$'s in the $\mathbf{LPC}\left(x,n\right)$ matrix is $n$, and the ${{n}^{2}}-n$ off-diagonal positions are occupied, in the equal number, by $x$ or $1/x$. Thus the Frobenius norm \eqref{eq:Frobenius} is: \\

\begin{equation}
\label{eq:LPC}
{{\left\| \mathbf{LPC}\left( x,n \right) \right\|}_{F}}={{\left( n+\frac{{{n}^{2}}-n}{2}\left( {{x}^{2}}+\frac{1}{{{x}^{2}}} \right) \right)}^{1/2}}.
\end{equation}
Despite the relation ${{\left( {Tr\left( {{\mathbf{A}}^{T}}\mathbf{A} \right)} \right)}^{1/2}}={{\left( n+\tfrac{{{n}^{2}}-n}{2}\left( {{x}^{2}}+\tfrac{1}{{{x}^{2}}} \right) \right)}^{1/2}}$ included in \eqref{eq:LPC}, there is no common formula for $CI$ of these matrices. As we will demonstrate, $\mathbf{LPC}\left( x,n \right)$ matrix with odd order $n$ has a maximal $CI$, which is false if the matrix order $n$ is an even number. 

\subsection{Case:	$n\ge 3$ odd}
\label{sec:LPCodd} 
For 	$n\ge 3$ odd, the column and row sums are: \\  \[\sum\nolimits_{i=1}^{n}{{{a}_{ij}}}=\sum\nolimits_{j=1}^{n}{{{a}_{ij}}}=1+\tfrac{n-1}{2}\left( x+\tfrac{1}{x} \right) \text{ for } i,j=1,...,n.\] 
It is easy to verify that the row sum is an eigenvalue with the corresponding eigenvector ${{\left( 1,1,...,1 \right)}^{T}}$. Considering that this eigenvalue equals the upper-bound mentioned in the above cited Theorem from \cite{AG1993}, this is the Perron root ${{\lambda}_{\max}}$ of $\mathbf{LPC}\left( x,n \right)$. The same result follows by Theorem 8.1.22. in \cite{HJ1985} since the row sums are equal. 
Consequently,  $\mathbf{LPC}\left(x,n\right)$ matrix for odd $n$ has maximal $CI$: 

\begin{equation}
\label{eq:CILPC_odd}
CI\left( \mathbf{LPC}\left( x,n \right) \right)=\frac{{{\left( x-1 \right)}^{2}}}{2x}.
\end{equation}
The above formula confirms the Lusk's guess, as states in \cite{AG1993}: ``Lusk \cite{Lusk1979} had further observed empirically that, at least for odd values of $n$, $\lambda_{max}$ seemed to reach a maximum at $1+(n-1)/2(S+1/S) \le  nS$.''

However, to the best of our knowledge, no one has noticed how important the existence and properties of matrix $\mathbf{LPC}$ are for normalization of $CI$.

It is important to notice that $CI\left(\mathbf{LPC}\left(x,n\right)\right)$ is independent of $n$ since the normalization does not depend on $n$. It is highly unfortunate that it was not known at the time of publishing \cite{Saaty1977}. \\

\noindent Let us examine the Gerschgorin disc and the location of eigenvalues.
When $n$ is odd, the matrix is circulant and its eigenvalues deduced from \eqref{eq:eval_circulant} are: \\
\[
{{\lambda}_{m}}=\sum\limits_{k=0}^{n-1}{{{c}_{k}}{{e}^{-2\pi imk/n}}=1+}\sum\limits_{k=1}^{n-1}{{{x}^{{{\left( -1 \right)}^{k+1}}}}{{e}^{-2\pi imk/n}}, m=1,...,n.}
\]
A straightforward computation leads to $n-1$ complex two by two conjugate eigenvalues with equal real parts $\operatorname{Re}{{\lambda }_{m}}=-\tfrac{{{\left( x-1 \right)}^{2}}}{2x}$, and one real positive ${{\lambda }_{\max }}=1+\tfrac{n-1}{2}\left( x+\tfrac{1}{x} \right)$. Thus, ${{\lambda}_{\max}}$ lies on the circle contouring the unique Gerschgorin disc 
\[D\left( {{a}_{ii}},{{R}_{i}} \right)=\left\{ z\in \mathbb{C}:\left| z-1 \right|\le \tfrac{n-1}{2}\left( x+\tfrac{1}{x} \right) \right\},i=1,...,n\]
of $\mathbf{LPC}\left( x,n \right)$, and all other eigenvalues are located inside the circle, on a~vertical line.

\begin{figure}[h]
	\centering
	\includegraphics[width=2.5in]{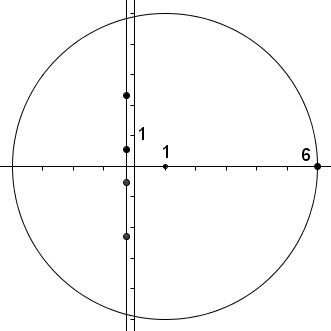}
	\caption[Case: $n=5$ and $x=2$]{Case: $LPC(x=2,n=5)$}
	\label{fig:gers52}
\end{figure}

Figure~\ref{fig:gers52} illustrates the situation for $n=5$ and $x=2$.

\subsection{Case:	$n\ge 4$ even}\label{LPCeven}
For $n\ge 4$ even, the row sums alternate. They give  $\sum\nolimits_{j=1}^{n}{{{a}_{ij}}}=1+\tfrac{n}{2}\left( x+\tfrac{1}{x} \right)-\tfrac{1}{x}$ for $i$ odd and  $\sum\nolimits_{j=1}^{n}{{{a}_{ij}}}=1+\tfrac{n}{2}\left( x+\tfrac{1}{x} \right)-x$ for $i$ even. It means that the matrix has two concentric Gerschgorin circles with radii: 
\begin{equation}
\label{eq:rR}
r(x,n)=\tfrac{n}{2}\left( x+\tfrac{1}{x} \right)-x \text{ and }  R(x,n)=\tfrac{n}{2}\left( x+\tfrac{1}{x} \right)-\tfrac{1}{x}.
\end{equation}
Assuming $x>1$ and with the use of Theorem 8.1.22. in~\cite{HJ1985} we get 
\begin{equation}
\label{eq:lambda_bounds}
1+r(x,n)\le {{\lambda }_{\max }}\le 1+R(x,n).
\end{equation}
By the Theorem quoted from \cite{AG1993}, $1+\tfrac{1}{2}\left( n-1 \right)\left( x+\tfrac{1}{x} \right)$ is an upper bound for ${{\lambda}_{\max}}$; we denote it $1+m(x,n)$. Notice that $m(x,n)$ is the arithmetic mean $A\left( r(x,n),R(x,n) \right)$ and that
\begin{equation}
\label{eq:diff_Rr}
R(x,n)-m(x,n)=m(x,n)-r(x,n)=\frac{{{x}^{2}}-1}{2x}.
\end{equation}
%A straightforward computation  shows
It is easy to verify (by simple computation) that $1+m(x,n)$ is not an eigenvalue, if $n$ is even. Therefore ${{\lambda }_{\max }}<1+m(x,n)$ and all eigenvalues are located inside Gerschgorin circle $\left| z-1 \right|=m(x,n)$. Starting from this point, with the use of \eqref{eq:CI},\eqref{eq:lambda_bounds} and \eqref{eq:diff_Rr}, we obtain:
\begin{equation}
\label{eq:CILPC}
0<\frac{{{(x-1)}^{2}}}{2x}-\frac{{{x}^{2}}-1}{2x(n-1)}\le CI\left( \mathbf{LPC}(x,n) \right)<\frac{{{(x-1)}^{2}}}{2x}.
\end{equation}

Inequality~\eqref{eq:CILPC} shows that the consistency index can be as near as we wish to its upper bound $\tfrac{{{(x-1)}^{2}}}{2x}$. This signifies that the large enough order $n$ of the matrix makes any error tolerable when $CI$ is used as inconsistency index and it is the case of \textit{reductio ad absurdum}.

As for ${{\lambda }_{\max }}$, numerical estimates collected in Table~\ref{table} give interesting information and inspire us to formulate our hypothesis~\ref{hypo}. 

\begin{table}[h]
\renewcommand{\arraystretch}{1.3}
\caption{$\mathbf{LPC}$ with even order: values of ${{\lambda }_{\max }}$ with its lower and upper bounds expressed as harmonic and geometric mean of $r$ and $R$}  \
\tabcolsep=2.5pt
\centering
\label{table}
\small
\begin{tabular}{|c||c|c|c||c|c|c||c|c|c|}
		\hline
&  & $n=4$ &  &  & $n=6$ &  & & $n=12$ &   \\ \hline
$x$ & $1+H$ & ${{\lambda}_{\max }}$ & $1+G$ & $1+H$ & ${{\lambda}_{\max }}$ & $1+G$ & $1+H$ & ${{\lambda}_{\max }}$ & $1+G$  \\ 
\hline
\hline
2 & 4.600000 & 4.644739 & 4.674235 & 7.160000 & 7.181572 & 7.204837 & 14.70909 & 14.71627 & 14.72953 \\ \hline
3 & 5.644444 & 5.762638 & 5.818944 & 9.120000 & 9.181326 & 9.225975 & 19.23636 & 19.25901 & 19.28478 \\ \hline
4 & 6.823529 & 7.015566 & 7.093029 & 11.29412 & 11.39625 & 11.45825 & 24.22459 & 24.26352 & 24.29968 \\ \hline
5 & 8.061538 & 8.324670 & 8.421590 & 13.55692 & 13.69842 & 13.77654 & 29.39860 & 29.45326 & 29.49912 \\ \hline
6 & 9.330330 & 9.662277 & 9.778129 & 15.86486 & 16.04441 & 16.13825 & 34.66585 & 34.73569 & 34.79102 \\ \hline
7 & 10.61714 & 11.01628 & 11.15090 & 18.19886 & 18.41549 & 18.52491 & 39.98649 & 40.07111 & 40.13582 \\ \hline
8 & 11.91538 & 12.38058 & 12.53392 & 20.54923 & 20.80228 & 20.92721 & 45.34056 & 45.43966 & 45.51369 \\ \hline
9 & 13.22132 & 13.75175 & 13.92380 & 22.91057 & 23.19954 & 23.33997 & 50.71693 & 50.83030 & 50.91363 \\ \hline
\end{tabular}
\end{table}

\begin{hypothesis}
\label{hypo}
The Perron root ${{\lambda }_{\max }}$ of a $\mathbf{LPC}\left( x,n \right)$ matrix with $n$ even is a function increasing in both variables $x$, $n$ and such that 
\begin{equation}
\label{eq:hyp}
1+H(r,R)\le {{\lambda }_{\max }}\le 1+G(r,R),
\end{equation}
where $r=r(x,n)$, $R=R(x,n)$ are defined in \eqref{eq:rR}, and $H(r,R)$, $G(r,R)$ are their harmonic and geometric means.
\end{hypothesis}

\section{Circulant PC matrix ($C_kPC$)}
%\label{}
Adding this section was motivated by the possibility of obtaining the exact value of Perron root and consequently a~close-form expression of consistency index. Consider reciprocal matrices that are circulant, as those studied in Section \ref{sec:LPCodd}.
First, we examine circulant matrices with the lowest inconsistency. Subsequently, we extend our results on other types of circulant reciprocal inconsistent matrices with entries $1$, $x$ and $1/x$.

\subsection{$C_kPC$ with low inconsistency}
Consider $n>3$, $x>0$, $x\ne 1$ and a circulant reciprocal $n\times n$ matrix with all terms ${{a}_{ij}}=1$ except one $x$ and one $1/x$ in each row. This determines the value $k=1$. Denote ${{\mathcal{F}}_{1}}$ this family of matrices and chose its representative ${{\mathbf{C}}_{\mathbf{1}}}\mathbf{PC}(x,n)=\left( {{a}_{ij}} \right)$ where 
$$
{{a}_{ij}}=\left\{ 
\begin{aligned}
& x\,\,\text{if}\,\,j-i=1\bmod n, \\ 
& \tfrac{1}{x}\,\,\text{if}\,\,i-j=1\bmod n, \\ 
& 1\,\,\text{otherwise}.
\end{aligned}
\right.
$$
For $n=5$, we have: 
\[{\mathbf{C_1}}\mathbf{PC}(x,5)=\begin{bmatrix}
1 & x & 1 & 1 & 1/x  \\
1/x & 1 & x & 1 & 1  \\
1 & 1/x & 1 & x & 1  \\
1 & 1 & 1/x & 1 & x  \\
x & 1 & 1 & 1/x & 1  \\
\end{bmatrix}.\]
Notice that sums of columns and rows of a circulant matrix are equal. For ${\mathbf{C_1}}\mathbf{PC}(x,n)$, the following sum is obtained: \[\sum\nolimits_{i=1}^{n}{{{a}_{ij}}}=\sum\nolimits_{j=1}^{n}{{{a}_{ij}}}=n-2+x+\tfrac{1}{x}=n+\tfrac{{{\left( x-1 \right)}^{2}}}{x} \text{ for } i,j=1,...,n.\] The row sum of a positive circulant matrix is its Perron root, as it follows from Theorem (already quoted) in \cite{HJ1985}. Hence ${{\lambda }_{\max }}=n+\tfrac{{{\left( x-1 \right)}^{2}}}{x}$, and the consistency index is:
\begin{equation}
\label{eq:CIC1PC}
CI\left({{\mathbf{C}}_{\mathbf{1}}}\mathbf{PC}(x,n)\right)=\frac{n+\tfrac{{{(x-1)}^{2}}}{x}-n}{n-1}=\frac{1}{n-1}\frac{{{(x-1)}^{2}}}{x}.
\end{equation}
The following two remarks correspond to the natural intuition: 
\begin{enumerate}
\item In a circulant matrix with large enough $n$, deviating entries (and their reciprocals) on a single diagonal does not have great impact on the inconsistency when comparing with a consistent $\mathbf{CC}$ matrix. Mathematical  formulation is:\\
%\item In a circulant matrix with large enough $n$,  deviating entries (and their reciprocals) on a one diagonal does not have great impact on the inconsistency when comparing with a consistent $\mathbf{CC}$ matrix. \\
\emph{For a fixed $x$, $\underset{n\to \infty }{\mathop{\lim }}\,CI\left( {{\mathbf{C}}_{\mathbf{1}}}\mathbf{PC}(x,n) \right)=\underset{n\to \infty }{\mathop{\lim }}\,\tfrac{1}{n-1}\tfrac{{{\left( x-1 \right)}^{2}}}{x}=0$.} 
\item PC matrix $\mathbf{CC}$ with all entries ${{a}_{ij}}=1$ is trivially consistent. A larger deviation of a PC matrix entry $x$ from $1$ results in a larger inconsistency (for a fixed $n$). Mathematically:\\
\emph{For a fixed $n$ and $x \in [1,\infty)$, ${{\mathbf{C}}_{\mathbf{1}}}\mathbf{PC}(x,n)$ is an increasing function in $x$.}
\end{enumerate}
It is worth to notice that $CI\left({\mathbf{C_1}}\mathbf{PC}(x,n)\right)=\tfrac{1}{n-1}\tfrac{{{\left( x-1 \right)}^{2}}}{x}$ hence it is vanishing as ${n\to \infty }$, and increasing for $x$ in $\left[ 1,\infty  \right)$.

The matrices of the family ${{\mathcal{F}}_{1}}$ differ only in the position of $x$ and $1/x$ in the rows hence the row sum is an invariant of ${{\mathcal{F}}_{1}}$. Consequently, the Perron root and its corresponding eigenvector ${{\left( 1,1,...,1 \right)}^{T}}$ (hence also $CI$) are invariants of ${{\mathcal{F}}_{1}}$. 

\subsection{$C_kPC$ with increased inconsistency}
The method used in the previous subsection can be extended to a family ${{\mathcal{F}}_{k}}$ of circulant reciprocal matrices with $k$ terms $x$ and equal number of $1/x$ in each row. Here $2k\le n-1$, because at least one $1$ (diagonal element) is obvious in each row.

As it was for ${{\mathcal{F}}_{1}}$, the family ${{\mathcal{F}}_{k}}$ has important invariants: the Perron root which is equal to the row sum, the corresponding eigenvector ${{\left( 1,1,...,1 \right)}^{T}}$, and the consistency index.
For $k\le \tfrac{n-1}{2}$, the unique row sum \[\sum\nolimits_{j=1}^{n}{{{a}_{ij}}}=n-2k+k\left( x+\tfrac{1}{x} \right)=n+k\tfrac{{{\left( x-1 \right)}^{2}}}{x}\] is necessarily a Perron root (see the above cited Theorem from \cite{HJ1985}). Consequently, $CI$ of each matrix in ${{\mathcal{F}}_{k}}$ is:
\begin{equation}
\label{eq:CIkPC}
CI\left( {{\mathbf{C}}_{k}}\mathbf{PC}(x,n) \right)=\frac{n+k\tfrac{{{(x-1)}^{2}}}{x}-n}{n-1}=\frac{k}{n-1}\frac{{{(x-1)}^{2}}}{x}.
\end{equation}
For a fixed $x$, a larger $k$ increases $CI$. The maximal inconsistency $\tfrac{{{\left( x-1 \right)}^{2}}}{2x}$ is reached for $k=\tfrac{n-1}{2}$ hence $n$ needs to be odd since $k$ is integer. We remind that this value of $CI$ was obtained in Section~\ref{sec:LPCodd} for $\mathbf{LPC}(x,n)$. It is not a~coincidence since $\mathbf{LPC}(x,n)$ matrices with odd order are members of ${{\mathcal{F}}_{\tfrac{n-1}{2}}}$.
 
If $n$ is even, the greatest value of $k$ is $\tfrac{n-2}{2}$. Thus has $CI\left( {{\mathbf{C}}_{k}}\mathbf{PC}(x,n) \right)$ the greatest possible value
\[CI\left( {{\mathbf{C}}_{\tfrac{n-2}{2}}}\mathbf{PC}(x,n) \right)=\tfrac{n-2}{n-1}\tfrac{{{(x-1)}^{2}}}{2x}.\] For a fixed $x$ and  $n\to \infty$, this value approaches the maximal inconsistency $\tfrac{{{\left( x-1 \right)}^{2}}}{2x}$.

\section{Full PC matrix (FPC)}
\label{sec:FPC}
Let us consider the following matrix $\mathbf{FPC}(x,n)=\left( {{a}_{ij}} \right)$:
$${{a}_{ij}}=\left\{ 
\begin{aligned}
& x\,\,\text{if}\,\,i<j, \\ 
& 1\,\,\text{if}\,\,i=j, \\ 
& 1/x\,\,\text{if}\,\,i>j,  
\end{aligned} 
\right.
$$ 
\noindent for $x>0$ and $n\ge 3$. Trivially,  $\mathbf{FPC}(x,n)$ is only consistent for $x=1$.

$\mathbf{FPC}(x,n)$ matrix is reciprocal with all entries above the main diagonal equal to $x$:
\[\mathbf{FPC}(x,n)=\begin{bmatrix}
1 & x & \cdots  & x & x  \\
1/x & 1 & \cdots  & x & x  \\
\vdots  & \vdots  & \ddots  & \vdots  & \vdots   \\
1/x & 1/x & \cdots  & 1 & x  \\
1/x & 1/x & \cdots  & 1/x & 1  \\
\end{bmatrix}.\]

Similarly to the $\mathbf{LPC}(x,n)$ matrix in Section~\ref{sec:LPC}, the number of $1$'s in $\mathbf{FPC}(x,n)$ is $n$, and the ${{n}^{2}}-n$ off-diagonal positions are occupied equally by $x$ or $1/x$.  Thus are their Frobenius norms equal:
\[{{\left\| \mathbf{FPC}\left( x,n \right) \right\|}_{F}}={{\left\| \mathbf{LPC}\left( x,n \right) \right\|}_{F}}={{\left( n+\frac{{{n}^{2}}-n}{2}\left( {{x}^{2}}+\frac{1}{{{x}^{2}}} \right) \right)}^{1/2}}.
\]
Recall from \eqref{eq:Frobenius} that the Frobenius norm is ${{\left\| \mathbf{A} \right\|}_{F}}={{\left({Tr\left( {{\mathbf{A}}^{T}}\mathbf{A} \right)} \right)}^{1/2}}$. Thus, one could expect the same eigenvalues for $\mathbf{FPC}(x,n)$ and $\mathbf{LPC}(x,n)$. However, unlike the circulant matrices, a different position of $x$ and $1/x$ influences the Perron root in $\mathbf{FPC}(x,n)$ matrix and thus affects the consistency index of the matrix.\\

\noindent Computations provided in \cite{KS2014} show that the Perron root of $\mathbf{FPC}(x,n)$ is equal to:
$$\lambda_n(x)= \frac{x-1}{x} \cdot \frac{x+x^{\frac{2}{n}}}{x^{\frac{2}{n}} -1},$$
and $$CI\left(\mathbf{FPC}(x,n)\right)=\frac{\lambda_n(x)-n}{n-1}=a_n(x)-\frac{n}{n-1},$$
where: 
$$a_n(x)=\frac{\lambda_n(x)}{n-1}=\frac{(x-1)(x^\frac{2}{n}+x)}{x}\cdot\frac{1}{n(x^\frac{2}{n}-1)-(x^\frac{2}{n}-1)}.$$
For $n \rightarrow \infty$, we have $$x^\frac{2}{n} \rightarrow 1$$ and $$n(x^\frac{2}{n}-1) \rightarrow 2\ln{x}$$ hence
$$a_n(x) \rightarrow \frac{x^2-1}{2x\ln{x}}$$
and 
\begin{equation}
\label{eq:limitCIFPC}
CI\left(\mathbf{FPC}(x,n)\right) \rightarrow \frac{x^2-1}{2x\ln{x}}-1.
\end{equation}
Evidently, for $x$ and $n$ large enough $CI\left(\mathbf{FPC}(x,n)\right)$ can be arbitrarily large, but it does not exceed nor reach the maximal value $\frac{{{\left( x-1 \right)}^{2}}}{2x}$ from \eqref{eq:CI_bounds}. 
However, for $x=2$ its rounded limit will be equal to $0.082$, hence will be AHP acceptable.
Moreover, when we consider small $n$, we can obtain the same result even for bigger $x$. \\

The following FPC matrix is taken from \cite{KS2014} and its eigenvalue-based inconsistency is computed as:
$$CI\left(\mathbf{FPC}(2.25,4)\right))=\frac{1}{18}\approx 0.055555556.$$
We can even increase $x$ to $3.375$ for $n=3$, which results in
$$CI\left(\mathbf{FPC}(3.375,3)\right))=\frac{1}{12}\approx 0.083333333,$$
which is still considered as acceptable by AHP theory. \\

Example~\ref{ex} brings a closed form of \ $CI\left(\mathbf{FPC}(x,3)\right)$ and shows that $CI\left(\mathbf{FPC}(x,3)\right)$ is smaller than $CI\left(\mathbf{LPC}(x,3)\right)$ for an arbitrary value of $x\neq1$. 
\begin{ex}
\label{ex}
\noindent Let $n\ge 3$, $x\ne 1$. The matrices \[\mathbf{FPC}(x,3)=\begin{bmatrix}
1 & x & x  \\
1/x & 1 & x  \\
1/x & 1/x & 1  \\
\end{bmatrix}  \text{and }  \mathbf{LPC}(x,3)=\begin{bmatrix}
1 & x & 1/x  \\
1/x & 1 & x  \\
x & 1/x & 1  \\
\end{bmatrix}\] differ in the position of unique entry $x$ and its reciprocal, which is the minimal difference.
Characteristic polynomials and Perron roots are \\
	$CharPol(\mathbf{FPC})=-{{\lambda }^{3}}+3{{\lambda }^{2}}+x+{{x}^{-1}}-2$, ${{\lambda }_{\max }}(\mathbf{FPC})={{x}^{\tfrac{1}{3}}}+{{x}^{-\tfrac{1}{3}}}+1,$ \\
	$CharPol(\mathbf{LPC})=-{{\lambda }^{3}}+3{{\lambda }^{2}}+{{x}^{3}}+{{x}^{-3}}-2$, ${{\lambda }_{\max }}(\mathbf{LPC})=x+{{x}^{-1}}+1$. \\
	
As we have seen in Section~\ref{sec:LPCodd}, 
$\mathbf{LPC}(x,n)$ achieves the maximal inconsistency 
$\tfrac{{{\left( x-1 \right)}^{2}}}{2x}$ for any odd $n\ge 3$. This is not the case of $\mathbf{FPC}(x,3)$, because 
\begin{equation}
\label{eq:CIFPC}
CI\left(\mathbf{FPC}(x,3)\right)=\frac{{{x}^{\tfrac{1}{3}}}+{{x}^{-\tfrac{1}{3}}}+1-3}{2}=\frac{{{\left( \sqrt[3]{x}-1 \right)}^{2}}}{2\sqrt[3]{x}}<\frac{{{\left( x-1 \right)}^{2}}}{2x}.
\end{equation}
The last inequality is a trivial consequence of the fact that the function $f(x)=\tfrac{{{\left( x-1 \right)}^{2}}}{2x}$ is decreasing on $\left( 0,1 \right]$ and increasing on $\left[ 1,\infty  \right)$. Thus, for any positive $x\ne 1$ is $f\left( \sqrt[3]{x} \right)<f\left( x \right)$, in our notation \[CI\left( \mathbf{FPC}(x,3) \right)<CI\left( \mathbf{LPC}(x,3) \right).\]
\end{ex}

\section{Conclusions}
In the present study, we have analyzed the Perron root and the consistency index $CI$ % \ref{eq:CI}
of special cases of PC matrices which are Toeplitz matrices with entries $1$, $x$ and $1/x$. The mathematical expressions for $CI$  \eqref{eq:CILPC_odd},\eqref{eq:CIC1PC},\eqref{eq:CIkPC},\eqref{eq:CIFPC} and their estimates in \eqref{eq:2odhadCICPC},\eqref{eq:CILPC} have values ranging from $0$ to the lowest upper bound $\tfrac{{{\left( x-1 \right)}^{2}}}{2x}$ inclusively. In the hypothesis~\ref{hypo}, a new estimate of $CI$ for a layer-cake PC matrix with even order is found. The presented results in this study are of considerable importance for PCs research.

It is worth noticing that the simplified version of pairwise comparisons,
presented in \cite{KS2015}, does not have inconsistency.
PC matrix elements are generated (from principal elements located on the diagonal
above or below the main diagonal) so that consistency condition is preserved.

\section*{Acknowledgments}
This research has been supported in part by the Euro Research grant ``Human Capital''. This paper was typeset using \LaTeXe. The authors acknowledge using the computational knowledge engine {WolframAlpha} \cite{WAlpha} and {wxMaxima} \cite{maxima} for the verification of derived results. The authors would also like to express appreciation to Tyler D. Jessup (Laurentian University, Computer Science), Grant O. Duncan, Team Lead, Business Intelligence, Integration and Development, Health Sciences North, Sudbury, Ontario, Canada), and Michal Zajac (Slovak University of Technology) for the editorial improvements.

\end{document}